\documentclass[12pt]{article}

\usepackage{amssymb}
\usepackage{amsthm}
\usepackage{bm}
\usepackage{amsmath}
\usepackage{graphicx}
\usepackage{cases}
\usepackage{amsfonts}

\usepackage[utf8]{inputenc}

\title{A stabilized local integral method using RBFs for the Helmholtz equation with applications to wave chaos and dielectric microresonators}

\author{L. Ponzellini Marinelli$^{[1,2]}$\\
        luciano@fceia.unr.edu.ar\\
        \and
        L. Raviola$^{[1]}$\\
        raviola@fceia.unr.edu.ar\vspace{.5cm}\\
        $[1]$ Faculty of Exact Sciences, Engineering and Surveying,\\
        National University of Rosario, Rosario, Argentina.\\
        $[2]$ Rosario Physics Institute, UNR-CONICET,\\
        Rosario, Argentina.\\
}
\date{\today}

\begin{document}
\maketitle

\begin{abstract}
Most problems in electrodynamics do not have an analytical solution so much effort has been put in the development of numerical schemes, such as the finite-difference method, volume element methods, boundary element methods, and related methods based on boundary integral equations. In this paper we introduce a local integral boundary domain method with a stable calculation based on Radial Basis Functions (RBF) approximations, in the context of wave chaos in acoustics and dielectric microresonator problems. RBFs have been gaining popularity recently for solving partial differential equations numerically, becoming an extremely effective tool for interpolation on scattered node sets in several dimensions with high-order accuracy and flexibility for nontrivial geometries. One key issue with infinitely smooth RBFs is the choice of a suitable value for the shape parameter which controls the flatness of the function. It is observed that best accuracy is often achieved when the shape parameter tends to zero. However, the system of discrete equations obtained from the interpolation matrices becomes ill-conditioned, which imposes severe limits to the attainable accuracy. A few numerical algorithms have been presented that are able to stably compute an interpolant, even in the increasingly flat basis function limit. We present the recently developed Stabilized Local Boundary Domain Integral Method in the context of boundary integral methods that improves the solution of the Helmholtz equation with RBFs. Numerical results for small shape parameters that stabilize the error are shown. Accuracy and comparison with other methods are also discussed for various case studies. Applications in wave chaos, acoustics and dielectric microresonators are discussed to showcase the virtues of the method, which is computationally efficient and well suited to the kind of geometries with arbitrary shape domains.
\end{abstract}


\section{Introduction and motivation}
\label{sec:intro}


Dielectric microresonators, also known as dielectric microcavities, have attracted interest in the last decades due to technological applications like microlasers and and as systems with intrinsic theoretical interest for its connections with quantum billiards and wave chaos \cite{cao_wiersig_2015,kaufman_kosztin_schulten_1999,wiersig2003}.

A quantum billiard is a system in which a free particle is confined within a 2D domain and whose dynamics is governed by the Schr\"odinger equation
\begin{equation}
 i\psi_t(\mathbf{x},t)=-\Delta \psi(\mathbf{x},t), \hspace{1cm} \mathbf{x}\in\Omega\subset\mathbb{R}^2, \;\; t>0.
 \label{schrodinger_eq}
\end{equation}
where $\psi(\mathbf{x},t)=0$ for $\mathbf{x}\in\Gamma$ being $\Gamma$ the boundary of the domain $\Omega$.

When searching the time harmonic solutions of this system in the form $\psi(\mathbf{x},t)=\tilde{\psi}(\mathbf{x})e^{ikt}$, the spatial dependence, $\tilde{\psi}(\mathbf{x})$, satisfies the well-known Helmholtz Equation (HE)
\begin{equation}
 \left( \Delta + k^2 \right)\tilde{\psi}(\mathbf{x})=0, \hspace{1cm} \mathbf{x}\in\Omega\subset\mathbb{R}^2, \;\; t>0.
 \label{helmholtz_eq}
\end{equation}
In this case, the eigenvalues to equation (\ref{helmholtz_eq}) are related to the energy of the particle.

On the other hand, a similar situation arises when trying to solve the problem of light waves propagating inside a dielectric medium satisfying the Maxwell equations. Also in this case, the search for time harmonic solutions leads to the Helmholtz equation for the spatial dependence of the electromagnetic field \cite{cao_wiersig_2015}.


For generic domains, the equation (\ref{helmholtz_eq}) cannot be solved analytically to find stationary states. So we must resort to finding efficient and reliable numerical methods to solve this equation. There are many numerical techniques to address this equation such as the finite element method (FEM), the finite volume method (FVM), the Boundary Element Method (BEM) or spectral methods (PS) \cite{trefethen_2000_book}. However, several of these require the construction of a specific mesh or refinement to efficiently address certain numerical problems on non-trivial geometries.

The BEM transforms the formulated Partial Differential Equations (PDE) into integral equations, that is, into an integral form over the boundary \cite{brebbia_dominguez_1998,partridge_brebbia_wrobel_1992}. In BEM the PDE that describes the physical problem is transformed into a Boundary Integral Equation (BIE), which is achieved by using Green's identities to then apply this integral formulation over points distributed in the domain. Many local integral methods are based on an integral formulation on small, strongly overlapping stencils with local interpolations.

In recent decades, methods involving the Radial Basis Functions (RBF) have become an extremely effective tool in non-trivial geometries for interpolation in sets of scattered nodes and for numerically approximating PDE. There are many modern books dealing with theory, implementations and applications \cite{fasshauer_2007_book,fasshauer_mccourt_2015_book,fornberg_flyer_2015_book}. One advantage is that when the distribution nodes are created, it is possible to achieve local refinement in critical areas depending on the specific problem \cite{fornberg_flyer_2015}. Particularly, this is interesting to resolve localized structures like the \emph{scarred} states observed in quantum chaos phenomena \cite{stockmann_1999_book}.

Using infinitely differential RBFs like Gaussians, exponential convergence can be shown. A practical obstacle is the ill-conditioning of the interpolation matrix when the shape parameter $\varepsilon$ that defines the Gaussian RBF tends to zero. It is known that when this parameter is reduced, the interpolation accuracy of the method improves considerably but the numerical conditioning of the problem worsens if it is solved with a direct type numerical method. That is, there is a conflict between accuracy and the constraint known as the uncertainty principle \cite{schaback_1995}.

In this paper we present the Stabilized Localized Boundary-Domain Integral Method (SLBDIM) \cite{ponzellini2021} in the context of Helmholtz type equations. This is a new stable integral local numerical method for approximating elliptic-type PDE solutions to solve Boundary Value Problems (BVP) in 2D that uses local interpolations with RBF for low values $\varepsilon >0$.
This technique is a combination of meshless methods, local integral formulations and boundary elements in multi-domains independent of a structured mesh and that only requires an unstructured distribution of nodes of the domain $\Omega$ and its boundary $\Gamma=\partial \Omega$ that allows to deal with complex geometries. For local interpolations, the Gaussian RBFs $\varphi (r)=e^{-(\varepsilon r)^2}$ are used when $\varepsilon\rightarrow 0$ in local interpolations in stable form.

Numerical results are shown for a small shape parameter that stabilizes the error. Comparisons with other methods in several cases are also discussed. It is shown that the method is computationally efficient and suitable for geometries that come from applications of wave chaos and dielectric microresonators. In particular, we solve differential problems with Dirichlet-type boundary conditions over square domains with quasi-uniform point distributions. 
 

\section{The Stabilized Localized Boundary Domain Integral Method for Helmholtz equations}
\label{sec:lbdim-st}

\subsection{Problem description and local integral method}

We consider the following Boundary Value Problem (BVP) on an open, bounded and simply connected domain $\Omega \subset \mathbb{R}^2$
 \begin{subnumcases}{(BVP)}
  \; \mathcal L \left[u\right](\mathbf{x}) = f(\mathbf{x}),  \hspace{1cm} \mathbf{x} \in \Omega, \label{ec:pvc_eliptico} \\
  \; \mathcal{B} \left[u\right](\mathbf{x}) = g(\mathbf{x}),  \hspace{1cm} \mathbf{x} \in \Gamma =\partial\Omega, \label{ec:pvc_eliptico_fro}
 \end{subnumcases}
where $\mathcal L[\ .\ ]=\Delta+\lambda $ is an elliptic differential Helmholtz-type operator, $\Delta = \frac{\partial}{\partial x^2}+\frac{\partial}{\partial y^2}$ is tha Laplacian, $\lambda \in \mathbb{R} $ (when $\lambda =k^2>0$, $k$ is the wave-number) and $f(\mathbf{x})$ is the smooth source term. $\mathcal B[\ .\ ]$ is the boundary operator with the boundary conditions (BC).

The BC are Dirichlet, Neumann or mixed over $\Gamma=\Gamma_1 \cup \Gamma_2$ and $\Gamma_1 \cap \Gamma_2 = \varnothing$
\begin{subnumcases}{}
 \hspace{.3cm} u(\mathbf{x}) = g_1(\mathbf{x}), \hspace{1.4cm} \mathbf{x}\in \Gamma_1, \label{ec:poisson_mixto_1}\\
 \hspace{0.15cm} \frac{\partial u(\mathbf{x})}{\partial n} = g_2(\mathbf{x}), \hspace{1.4cm} \mathbf{x}\in \Gamma_2, \label{ec:poisson_mixto_2x}
\end{subnumcases}
with $g_1$ and $g_2$ known data, and $\frac{\partial u(\mathbf{x})}{\partial n}$ the outward normal derivative of the unknown field $u$.

We propose that PDE (\ref{ec:pvc_eliptico}) can be written as
\begin{equation}
\Delta u \left( \mathbf{x}\right) = f(\mathbf{x})-\lambda u \left( \mathbf{x}\right) = 
b\left(\mathbf{x},u\left(\mathbf{x}\right)\right),
\label{ec:gov_equation}
\end{equation}
where $u \left(\mathbf{x}\right)$ is the potential in the point $\mathbf{x} \in \Omega$.

We consider $\mathbf{x}\in\Omega\subset\mathbb{R}^2$
\begin{equation}
\Delta u^*= \delta(\mathbf{x-{\mathbf \xi}}),
\label{ec:delta_potential}
\end{equation}
where $\delta(\mathbf{x-{\mathbf \xi}})$ is Delta's delta centered at $\mathbf \xi\in\Omega$ with fundamental solution
\begin{equation}
 u^*(\mathbf{x,{\mathbf \xi}}) = \frac{1}{2\pi}ln(r), \;\;\;\; r=\| \mathbf{x-{\mathbf \xi}}\|.
 \label{ec:sol_funda_2d}
\end{equation}

From equation (\ref{ec:gov_equation}) 
\begin{equation}
 \Delta u \left( \mathbf{x}\right) = 
 b \Leftrightarrow \int_{\Omega}u^{*}\left(\mathbf{x},{\mathbf \xi}\right)\Delta u\left(\mathbf{x}\right) \ d\Omega_{\mathbf{x}} = \int_{\Omega}u^{*}(\mathbf{x},{\mathbf \xi}) b \ d\Omega_{\mathbf{x}}.
 \end{equation}

Applying Green's second identity for $u$ that satisfies (\ref{ec:gov_equation}) and $u^*$ that satisfies (\ref{ec:delta_potential}) 
  \begin{equation}
\int_{\Omega} \left(u^{*} \Delta u - u\Delta u^{*} \right) d\Omega_{\mathbf{x}} = \oint_{\Gamma} \left( u^{*} \frac{\partial u}{\partial n} - u\frac{\partial u^*}{\partial n} \right) d\Gamma_{\mathbf{x}},
\label{ec:poisson_drm_2}
\end{equation}
we obtain 
\begin{equation}
 u({\mathbf \xi}) = \int_{\Omega}u^{*}\left(\mathbf{x},{\mathbf \xi}\right)b \ d\Omega_{\mathbf{x}} - \oint_{\Gamma}\left[ u^{*}\left(\mathbf{x},{\mathbf \xi}\right) \frac{\partial u(\mathbf{x})}{\partial n} - u(\mathbf{x})\frac{u^*(\mathbf{x},{\mathbf \xi})}{\partial n} \right] d\Gamma_{\mathbf{x}}.
\label{ec:ec_int_front_poisson_0}
\end{equation}

From equation (\ref{ec:ec_int_front_poisson_0}) we have a formula for the integral representation of the PDE over a subregion $\Omega_i$ with boundary $\Gamma_i$. The interior collocation point ${\mathbf \ xi}$ is obtained as before from the fundamental solution and Green's second identity
\begin{equation}
u({\mathbf \xi}) = 
\int_{\Gamma_i}q^{*}\left(\mathbf{x},{\mathbf \xi}\right)u\left(\mathbf{x}\right) \ d\Gamma_{\mathbf{x}} -
\int_{\Gamma_i}u^{*}\left(\mathbf{x},{\mathbf \xi}\right)q\left(\mathbf{x}\right) \ d\Gamma_{\mathbf{x}} +
\int_{\Omega_i} b \ u^{*}\left(\mathbf{x},{\mathbf \xi}\right)\  d\Omega_{\mathbf{x}},
\label{ec:int_eq_conv-diff}
\end{equation}
where $q=\frac{\partial u}{\partial n}$ is the normal derivative of the unknown field, $u^{*}$ is the fundamental Laplacian solution and $q^{*}= \frac{\partial u^{*}}{\partial n}$ is the normal derivative of the fundamental solution.

Using the well-known Green-Dirichlet function (FGD), $G\left(\mathbf{x},{\mathbf \xi}\right)$, and its normal derivative $Q\left(\mathbf{ x},{\mathbf \xi}\right)$ \cite{greenberg_2015} in (\ref{ec:int_eq_conv-diff}) we obtain a new integral formulation of the form
\begin{equation}
u({\mathbf\xi}) = \int_{\Gamma_i} Q\left(\mathbf{x},{\mathbf \xi}\right)u\left(\mathbf{x}\right) \ d\Gamma_{\mathbf{x}} + \int_{\Omega_i} b \ G\left(\mathbf{x},{\mathbf \xi}\right) \ d\Omega_{\mathbf{x}}.
\label{ec:int_eq_conv-diff(Green)}
\end{equation}
since the integral over $\Gamma_i$ involving $G$ in (\ref{ec:int_eq_conv-diff}) vanishes since its value is zero.

In addition, if the non-homogeneous term $b$ of the PDE can be split
\begin{equation}
b\left(\mathbf{x},u\left(\mathbf{x}\right)\right) = f\left(\mathbf{x}\right) - \lambda u\left(\mathbf{x}\right), 
\label{Eq:gov_equation_split}
\end{equation}
where the funcion source $f$ is data.

The integral representation (\ref{ec:int_eq_conv-diff(Green)}) in each subregion of integration $\Omega_i$ is
\begin{equation}
u({\mathbf\xi}) = \int_{\Gamma_i} Q(\mathbf{x},{\mathbf\xi}) u(\mathbf{x}) \ d\Gamma_{\mathbf{x}} + \int_{\Omega_i}  G(\mathbf{x},{\mathbf\xi}) f(\mathbf{x}) \ d\Omega_{\mathbf{x}} + \int_{\Omega_i} - \lambda u\left(\mathbf{x}\right) \ G(\mathbf{x},{\mathbf\xi}) \ d\Omega_{\mathbf{x}},
\label{Eq:int_eq_conv-diff(Green)_milfd}
\end{equation}
where ${\mathbf\xi}$ is the interior source point. The collocation technique is done only at interior points of the domain.

\subsection{Local interpolations with RBF}

A function $\varphi:\mathbb{R}^d\rightarrow \mathbb{R}$ is an RBF if there exists $\phi:[0,\infty)\rightarrow \mathbb{R}$ such that
\begin{equation}
\varphi\left(\mathbf{x}\right) = \phi(r), \hspace{1cm} r=\|\mathbf{x}-\mathbf{x}_j\|,
\label{ec:def_rbf_0}
\end{equation}
where $\| . \|$ is the Euclidean norm on $\mathbb{R}^d$ and depends on the distance to a center $\mathbf{x}_j\in\mathbb{R}^d$. If it depends on the shape parameter $\varepsilon>0$, then $\varphi_j^\varepsilon\left(\mathbf{x}\right) = \phi(r,\varepsilon)$ is often noted.

In the LBDIM the field $u$ is locally interpolated with RBF $\{\varphi_j\}_{j=1}^n$ with centers of the stencil $\Theta_{\mathbf{x}}=\{\mathbf{x}_j\}_{j=1}^n$
\begin{equation}
u\left(\mathbf{x}\right) \approx \sum^{n}_{j=1}\alpha_j\varphi_j(\mathbf{x}),
 \label{ec:interpola_u_milfd}
\end{equation}
where the interpolation matrix $\mathbf{A}_i$ is such that
\begin{equation}
 \left(\mathbf{A}_i\right)_{jk}=\varphi_k(\mathbf{x}_j)=\phi(\|\mathbf{x}_j-\mathbf{x}_k\|), \;\;\;j,k=1,\dots,n
\end{equation}

The term $b$ of (\ref{Eq:gov_equation_split}) is interpolated with RBF $\{\chi_j\}_{j=1}^m$ with centers of the stencil $\Theta_{\mathbf{y }}=\{\mathbf{y}_j\}_{j=1}^m$
\begin{equation}
\widetilde{b}\left(u\left(\mathbf{x}\right),\nabla u \left(\mathbf{x}\right) \right)\approx \sum_{j=1}^{m} \beta_j \chi_j\left(\mathbf{x}\right),
\label{ec:int_term_no_hom_milfd}
\end{equation}
where the interpolation matrix $\mathbf{\widetilde{A}}_i$ is such that
\begin{equation}
(\mathbf{\widetilde{A}}_i)_{jk}=\chi_k(\mathbf{y}_j) =\chi(\|\mathbf{y}_j-\mathbf{y}_k\|), j,k=1,\dots,m
\end{equation}
The RBFs are eventually of the same type and with the same centers. If we take the same RBF bases with the same centers, the result is $\{\varphi_j\}_{j=1}^n$ and $\{\chi_j\}_{j=1}^m$ for $m=n$ although they could be different depending on the application problem or numerical experience.

The local integral formulation of (\ref{Eq:int_eq_conv-diff(Green)_milfd}) is of the form
\begin{eqnarray}
u({\mathbf\xi}) &\approx& \sum_{j=1}^{n} \alpha_j \left\{\int_{\Gamma_i} Q(\mathbf{x},{\mathbf\xi}) \varphi_j(\mathbf{x}) \ d\Gamma_{\mathbf{x}} \right\} \nonumber \\ 
&+&\sum_{j=1}^{m} \beta_j \left\{\int_{\Omega_i} G(\mathbf{x},{\mathbf\xi})\chi_j\left(\mathbf{x}\right) d\Omega_{\mathbf{x}}\right\} + \int_{\Omega_i} G(\mathbf{x},{\mathbf\xi}) f(\mathbf{x}) \ d\Omega_{\mathbf{x}}.
\label{ec:DGF-DRM-meshless}
\end{eqnarray}

If $\Theta=\{ \mathbf{x}_1,\dots ,\mathbf{x}_N \}$ is the discretization of domain $\Omega$ and ${\mathbf\xi}=\mathbf{x} _i\in\Theta$ is the collocation point, the discretized formulae of the unknown field is
\begin{equation}
u_i=u\left(\mathbf{x}_i\right)= \sum^n_{j=1}\alpha_j \widetilde{h}_{ij} + \sum^m_{j=1}\beta_j \widetilde{g}_{ij} + \widetilde{f}_{i},
 \label{ec:disc_form_milfd}
\end{equation}
where $\alpha_j$ and $\beta_j$ come from equations (\ref{ec:interpola_u_milfd}) and (\ref{ec:int_term_no_hom_milfd}). The coefficients $\widetilde{h}_{ij}$, $\widetilde{g}_{ij}$ and $\widetilde{f}_{i}$ are of the form
\begin{subnumcases}{}
\widetilde{h}_{ij} = \int_{\Gamma_i} Q\left(\mathbf{x},\mathbf{x}_i\right)
\varphi_j\left(\mathbf{x}\right)d\Gamma_{\mathbf{x}}, \label{local_int_coeff_milfd_hij} \\
\widetilde{g}_{ij} = \int_{\Omega_i}G\left(\mathbf{x},\mathbf{x}_i\right)
\chi_j\left(\mathbf{x}\right)d\Omega_{\mathbf{x}}, \label{local_int_coeff_milfd_gij} \\
\widetilde{f}_i    = \int_{\Omega_i}G\left(\mathbf{x},\mathbf{x}_i\right)
f\left(\mathbf{x}\right)d\Omega_{\mathbf{x}},
\label{local_int_coeff_milfd_fi}    
\end{subnumcases}
which are calculated by Gauss-Legendre quadratures.

Defining the vectors ${\mathbf \alpha}=\left[\alpha_1,\dots,\alpha_n\right]^T$ and ${\mathbf \beta}=\left[\beta_1,\dots,\beta_m \right]^ T$ as interpolation coefficients, the discretized form (\ref{ec:disc_form_milfd}) of $u$ can be expressed as
\begin{equation}
   u_i = \mathbf{\widetilde{h}}_i^T {\mathbf \alpha} + \mathbf{\widetilde{g}}_i^T {\mathbf \beta} + \widetilde{f}_i,
   \label{ec:discret_form_ui_milfd_matrix_form}
\end{equation}
where $\mathbf{\widetilde{h}}_i=[\widetilde{h}_{i1},\dots,\widetilde{h}_{in}]^T$ and $\mathbf{\widetilde{g} }_i=[\widetilde{g}_{i1},\dots,\widetilde{g}_{im}]^T$ are the influence coefficients, and $\widetilde{f}_i\in\mathbb{R}$ is data.

The vector ${\mathbf \alpha}$ arises from the local system by interpolating with the RBF basis $\{\varphi_j\}_{j=1}^n$
\begin{equation}
 \mathbf{A}_i {\mathbf \alpha} = {\mathbf{d}}_i \Leftrightarrow {\mathbf \alpha} = \mathbf{A}_i^{-1}\mathbf{d}_i
\label{ec:int_u_milfd_sel_alfa}
\end{equation}
and the vector ${\mathbf \beta}$ arises from the local system by interpolating with the RBF basis $\{\chi_j\}_{j=1}^m$
\begin{equation}
 \mathbf{\widetilde{A}}_i {\mathbf \beta} = \mathbf{\widetilde{b}}_i \Leftrightarrow  {\mathbf \beta} = \mathbf{\widetilde{A}}_i^{-1}\mathbf{\widetilde{b}}_i = \mathbf{\widetilde{A}}^{-1}_i\left(\mathbf{A}_{\widetilde{b}_i}{\mathbf\alpha}\right) = \mathbf{\widetilde{A}}^{-1}_i\left(\mathbf{A}_{\widetilde{b}_i}\mathbf{A}^{-1}_i \mathbf{d}_i\right),
\label{ec:int_term_no_hom_milfd_sel}
\end{equation}
where $\mathbf{A}_{\widetilde{b}_i}$ is the calculation matrix of the vector $\mathbf{\widetilde{b}}_i$ with known coefficients 
\begin{equation}
(\mathbf{A}_{\widetilde{b}_i})_{jk}=\widetilde{b}\left(\varphi_k\left(\mathbf{y}_j\right),\nabla\varphi_k\left(\mathbf{y}_j\right)\right), \;\;\;\; j=1,\dots,m, k=1,\dots,n.
\label{ec:calc_matrix_term_no_hom}
\end{equation}

Substituting (\ref{ec:int_u_milfd_sel_alfa}) and (\ref{ec:int_term_no_hom_milfd_sel}) in the discretized form (\ref{ec:discret_form_ui_milfd_matrix_form}), we obtain the discretized matrix form for $u_i$ in terms of $\mathbf{d}_i$
\begin{equation}
u_i =  \left(\mathbf{\widetilde{h}}_i^T \mathbf{A}^{-1}_i +
 \mathbf{\widetilde{g}}_i^T \mathbf{\widetilde{A}}^{-1}_i \mathbf{A}_{\widetilde{b}_i}\mathbf{A}^{-1}_i\right)\mathbf{d}_i + \tilde{f}_i.
\label{local_int_eq-ii_milfd}
\end{equation}

Rewriting (\ref{local_int_eq-ii_milfd}) we obtain an algorithmic procedure to avoid the computation of inverses $\mathbf{A}_i^{-1}$ and $\mathbf{\widetilde{A}} _i^{-1}$ (see \cite{ponzellini2021_tesis})
\begin{equation}
 u_i = \mathbf{z}^T\mathbf{d}_i + \widetilde{f}_i \hspace{1cm}\mbox{donde} \;\; \mathbf{z}^T = \mathbf{\widetilde{h}}_i^T \mathbf{A}_i^{-1} + \mathbf{\widetilde{g}}_i^T \mathbf{\widetilde{A}}_i^{-1} \mathbf{A}_{\widetilde{b}_i} \mathbf{A}_i^{-1}
\label{local_int_eq-ii_zi}
\end{equation}
which are assembled into a global sparse-like system and numerically resolved with Generalized Minimal Residual (GMRES). 



\subsection{Stability with Gaussian RBFs}
\label{subsec:stability}

Convergence in global interpolations with $\varepsilon$-dependent RBFs can be studied in a stationary way ($n=cte.$ and $\varepsilon\rightarrow 0$) or non-stationary ($\varepsilon=cte.$ and $ increases n$). In the case of Gaussian RBFs, they produce convergence of order $O(e^{-\frac{const}{(\varepsilon h)^2}})$ (superspectral).


The RBF interpolation matrix is
\begin{equation*}
\mathbf{A}(\varepsilon) = 
\left[
\begin{array}{cccc}
\phi(\|\mathbf{x}_1-\mathbf{x}_1\|,\varepsilon) & \phi(\|\mathbf{x}_1-\mathbf{x}_2\|,\varepsilon) & \dots & \phi(\|\mathbf{x}_1-\mathbf{x}_{n}\|,\varepsilon)\\
\phi(\|\mathbf{x}_2-\mathbf{x}_1\|,\varepsilon) & \phi(\|\mathbf{x}_2-\mathbf{x}_2\|,\varepsilon) & \dots & \phi(\|\mathbf{x}_2-\mathbf{x}_{n}\|,\varepsilon)\\
\vdots & \vdots & \ddots & \vdots \\
\phi(\|\mathbf{x}_{n}-\mathbf{x}_1\|,\varepsilon) & \phi(\|\mathbf{x}_{n}-\mathbf{x}_2\|,\varepsilon) & \dots & \phi(\|\mathbf{x}_{n}-\mathbf{x}_{n}\|,\varepsilon)
\end{array}
\right].
\label{interpolation_matrix_eps}
\end{equation*}

When $\varepsilon$ is small, the RBFs become almost linearly dependent ('flat') forming a bad basis of functions and generating ill-conditioned interpolation matrices $\mathbf{A}(\varepsilon)$ in a good interpolation space. To avoid this problem in \cite{fornberg_larsson_flyer_2011,larsson_lehto_heryudono_fornberg_2013} numerical techniques were developed that stabilize the solutions of linear systems where the RBFs that form the matrix of the system take arbitrarily small shape parameters. The RBF-QR method developed for global interpolations of scattered nodes using Gaussian RBFs is numerically stable for nearly zero parameters. The idea of the RBF-QR algorithm is to change the base $\{\phi_j\}$ to a new base $\{\psi_j\}$ using combinations of polynomial powers, Chebyshev polynomials and trigonometric functions.

\section{Implementation of the SLBDIM}

The new matrix form for $u$ of (\ref{local_int_eq-ii_milfd}) at each node is
\begin{equation}
u_i = \left(\boldsymbol{l}_i^T {\boldsymbol{B}_i}^{-1} +
 \boldsymbol{\widetilde{l}}_i^T {\boldsymbol{{\widetilde{B}}}_i}^{-1} {\boldsymbol{{B}}_{\tilde{b}_i}}
 {\boldsymbol{{B}}_i}^{-1}\right)\boldsymbol{d}_i + \tilde{f}_i,
\label{local_int_eq-ii_rbf-qr}
\end{equation}
where $\boldsymbol{l}_i=\left[\ldots,l_{ik},\ldots\right]^T$ and $\boldsymbol{\widetilde{l}}_i=[\ldots,\widetilde{l}_{ik},\ldots]^T$ are the column vectors.

For internal stencils, the local interpolation matrix is 
\begin{equation}
{\boldsymbol{{B}}_\psi^i} = \boldsymbol{V}
 \left[\begin{array}{c}
 \boldsymbol{I}_n \\
 \boldsymbol{\widetilde{R}}^T \\
 \end{array}
 \right],
 \label{rbf-qr_collocation_matrix}
\end{equation}
where $({\boldsymbol{{B}}_\psi^i})_{jk}={\psi_k(\boldsymbol{x}_j)}$ and $V_{jk}=V_k(\boldsymbol{x}_j)$ for $j,k=1,\dots,n$ (\cite{fornberg_larsson_flyer_2011} for details).

For boundary stencils, the local matrix interpolation matrix is $\boldsymbol{B}^i$ has two blocks,
\begin{equation}
{ \boldsymbol{{B}}_i} = 
 \left[\begin{array}{c}
 {\boldsymbol{{B}}_\psi^i}\\
 {\boldsymbol{{B}}_{\mathcal B\psi}^i}\\
 \end{array}
 \right],
 \label{rbf-qr_matrix}
\end{equation}
where the first matrix block is
\begin{equation}
({\boldsymbol{{B}}_\psi^i})_{jk}={\psi_k(\boldsymbol{x}_j)}, \end{equation}
for $j=1,\dots,n_{int}$ (interior nodes) and $k=1,\dots,n$ (boundary nodes), and the second matrix block is
\begin{equation}
({\boldsymbol{{B}}_{\mathcal B\psi}^i})_{jk}=\mathcal B {\psi_k(\boldsymbol{x}_j)}
\end{equation}
for $j=n_{int}+1,\dots,n$ and $k=1,\dots,n$.

To avoid calculating ${\boldsymbol{{B}}_i}^{-1}$ and ${\boldsymbol{{\widetilde{B}}}_i}^{-1}$ when $\varepsilon\rightarrow 0$ we follow an algorithmic procedure. The inclusion of this technique in the local integral method allows to stabilize the numerical error of the approximation of the Helmholtz-type equations. This Stabilized Domain and Boundary Local Integral Method (SLBDIM) was presented at \cite{ponzellini2021} for Poisson problems, convection-diffusion equations and elliptic PDEs. Another strategy of stability technique for local integral methods that uses RBF interpolation functions was presented in \cite{ponzellini2021_revuma}.

\section{Numerical examples on several billiars}
\label{sec:billiars}

In this section we report two numerical experiments to show the accuracy and efficiency of the proposed numerical scheme to solve Helmholtz-type equations in two dimensions. Implementations and numerical experiments were performed using MATLAB version R2017a numerical calculation software on a PC with 7.5 GB of RAM and an Intel Core i7-7500U 7th Generation CPU. running at 2.70GHz. 

The reported errors are the standard error $L_2$ ($L_2$-Error)
\begin{equation}
\begin{array}{rcl}
 L_2\mbox{-Error} &=& \sqrt{\frac{\sum_{i=1}^{N}\left(u^{exac}_i-u^{approx}_i\right)^2}{ \sum_{i=1}^{N}\left(u^{exac}_i\right)^2}}
 \label{Eq:L2error}
\end{array}
\end{equation}
 and the root mean square error (RMS):
\begin{equation}
\begin{array}{rcl}
\mbox{RMS} &=& \sqrt{\frac{\sum_{i=1}^{N}\left(u^{exac}_i-u^{approx}_i\right)^2}{N} }.
\label{Eq:RMS}
\end{array}
\end{equation}

\subsection{Polygonal billiars: case 1}
\label{subsec:polygo1}

This Helmholtz-type PDE is given over the rectangular domain $\Omega = [-1, 1] \times [-1,1]$
\begin{equation}
\left\{
\begin{array}{rcl}
\Delta u(\boldsymbol{x})-k^2u(\boldsymbol{x}) &=& f(\boldsymbol{x}), \hspace{.5cm}  \boldsymbol{x}=(x,y)\in \Omega,\\
u(\boldsymbol{x}) &=& g(\boldsymbol{x}), \hspace{.5cm} (x,y)\in \Gamma=\partial\Omega,
\end{array}
\right. 
\label{Eq:edp_helmohltz_1}
\end{equation}
where $f(x,y)=2\cos (x^2+y)-(4x^2+1+k^2)sin(x^2+y)$ and the parameter $k=9 $. The BCs of this BVP are of the Dirichlet type, the analytical solution being $u(x,y)=sin(x^2+y)$. In our case, we will use the local integral method presented in its original form with Gaussian RBF kernels $\phi(r)=e^{-(\varepsilon r)^2}$ (we will call it LBDIM) and in its stabilized form (SLBDIM).

There are several ways to discretize the $\Omega$ domain with distributions of nodes. In our case we will use the algorithm for generating quasi-uniform distributions developed in \cite{fornberg_flyer_2015} for 2D. These distributions were created with a fast-forward method that generates a set of nodes from a density function starting from the $\Gamma$ boundary towards the interior of the domain.

\begin{figure*}[htb]
\includegraphics[scale=.45]{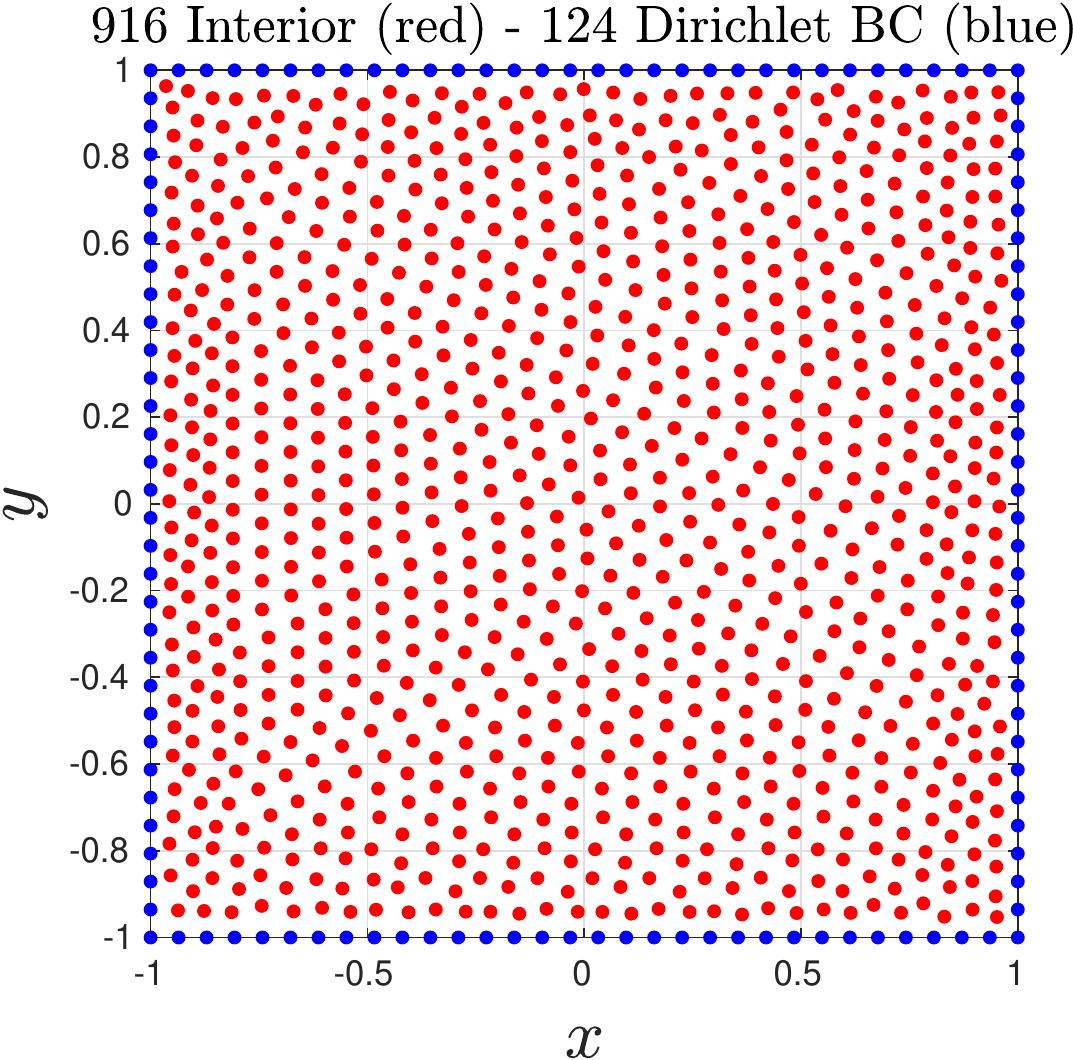}
\hspace{1cm}
\includegraphics[scale=.5]{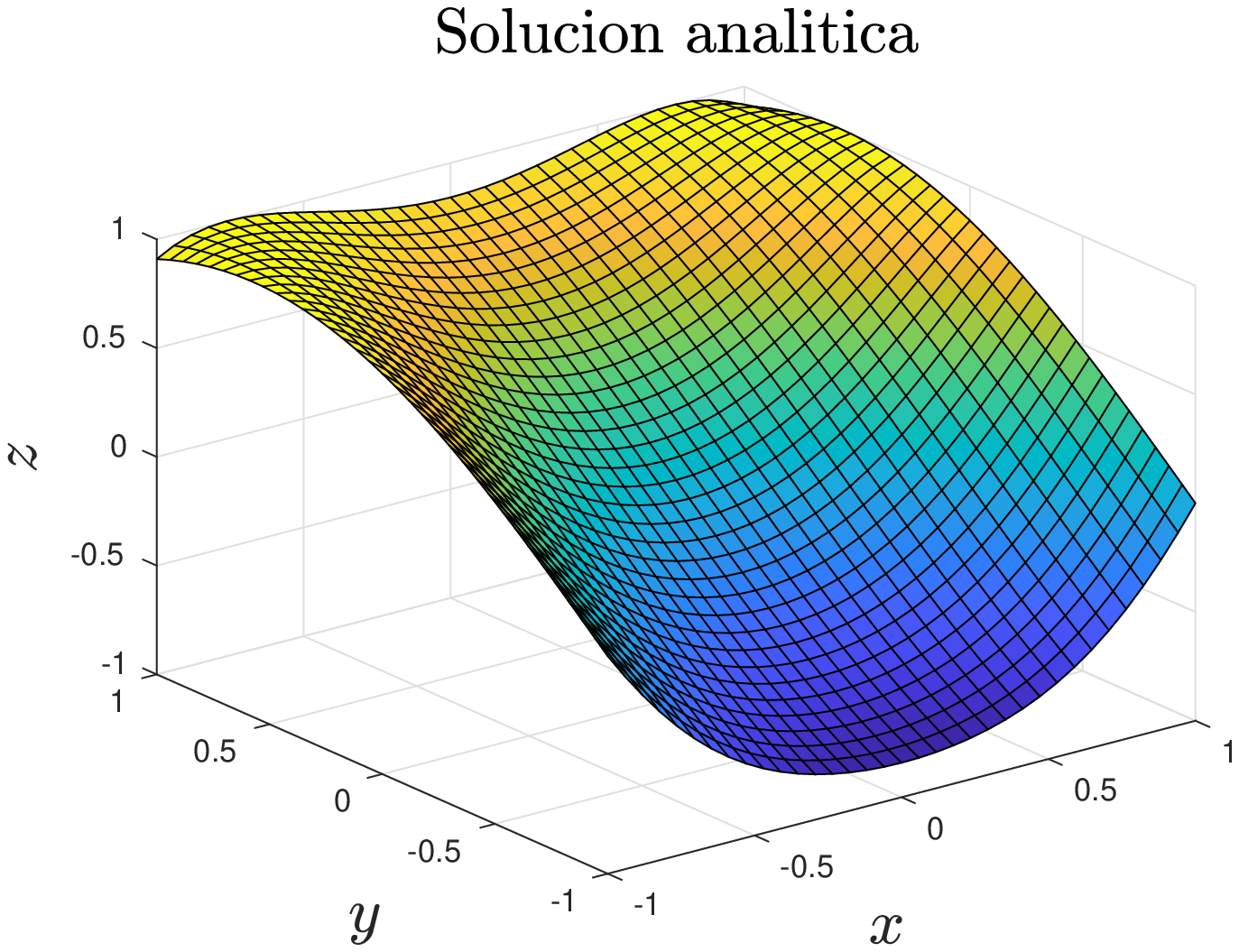}
\caption{Quasi-uniform 2D node distribution for $N_{int}=916$ internal collocation points and $N_{col}=124$ boundary points with Dirichlet BC (left). Analytical solution of BVP (right).}
\label{fg:HE_P7_quni400_quni916}
\end{figure*}

\begin{figure}[htb]
\centerline{\includegraphics[scale=.75]{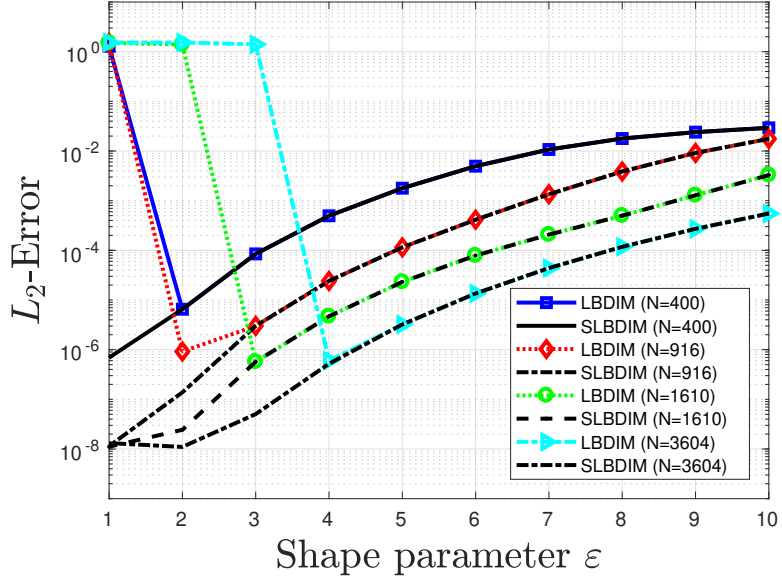}} \caption{Comparison of the $L_2$-Error between LBDIM and SLBDIM versus the shape parameter $\varepsilon$.}
\label{fg:HE_P7_ecm_vs_sp}
\end{figure}

We compare the $L_2$-Error of the formulation of the LBDIM and the SLBDIM using the Gaussian RBFs in the local interpolations varying the parameter in the form $\varepsilon \in [1,10]$. Figure \ref{fg:HE_P7_ecm_vs_sp} shows that as $\varepsilon$ decreases, the accuracy increases but the LBDIM is destabilized and the convergence is interrupted all for cases $N=400,916,1610,3604$ quasi-uniform nodes. In turn, we observe that as we increase the number of nodes on the domain and the boundary, the errors decrease. This plot shows that for local interpolation with Gaussian RBF lead to a loss in accuracy for small shape parameters. However, the best performance is obtained by the stabilized local integral method to address this Helmholtz-type equation with known analytical solutions. The error for $N=916,1610,3604$ is of order $1\times 10^{-8}$. The application of the RBF-QR kernel makes the system well-posed to solve them with a direct method in the LBDIM. In this numerical experiment the size of the stencil is $n=50$.

\begin{figure*}[htb!!!!]
\includegraphics[scale=.45]{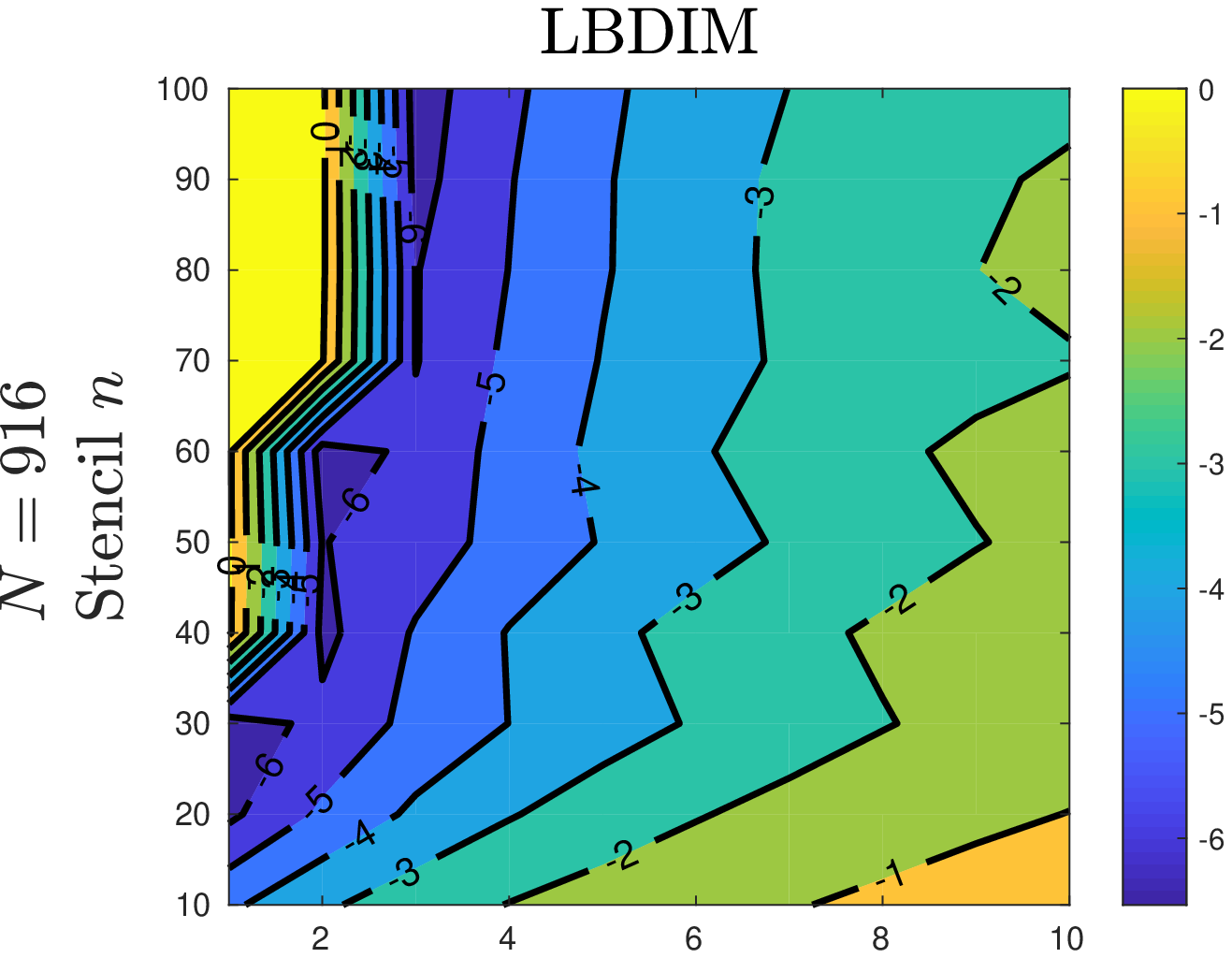}
\includegraphics[scale=.45]{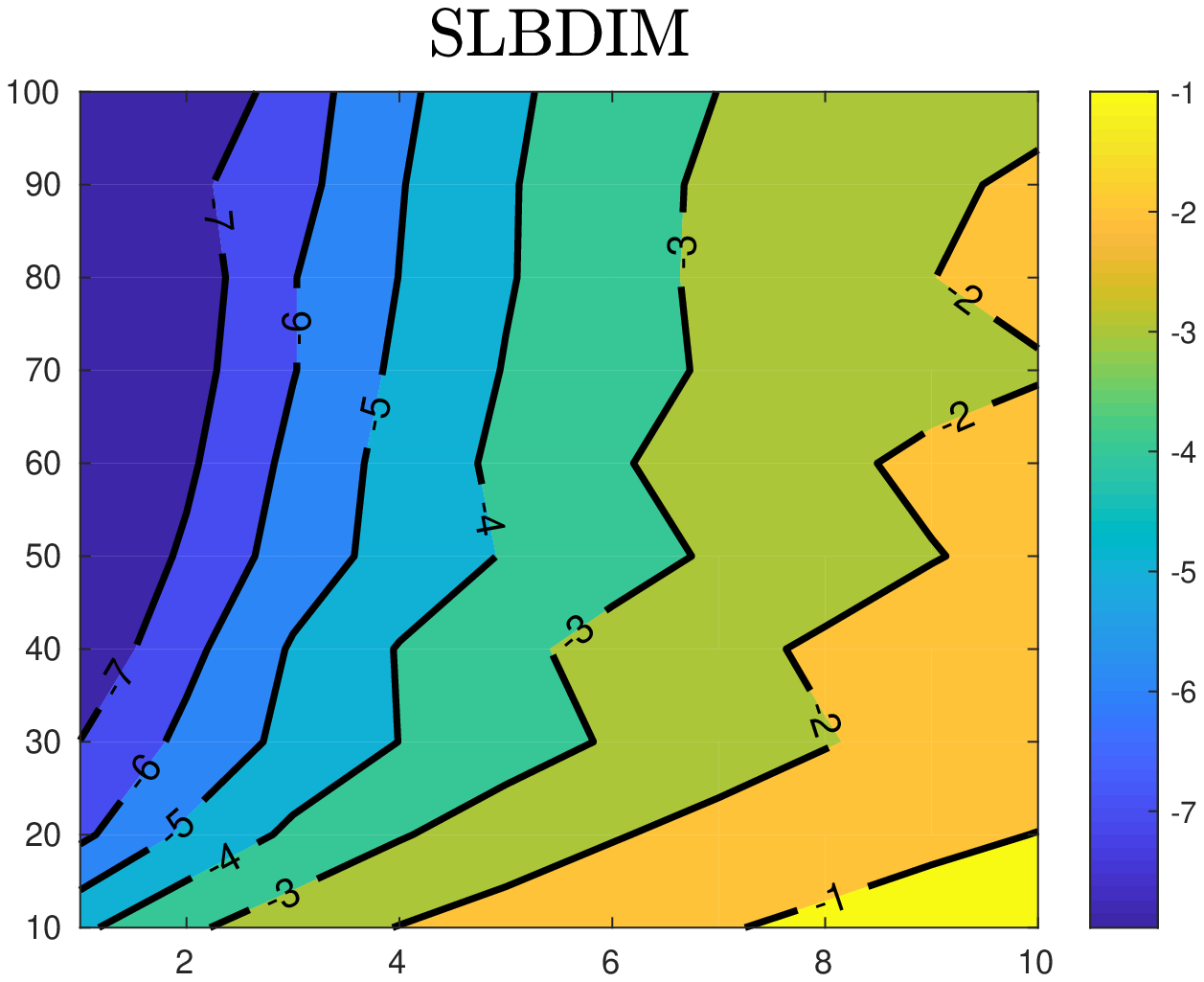}\\
\includegraphics[scale=.45]{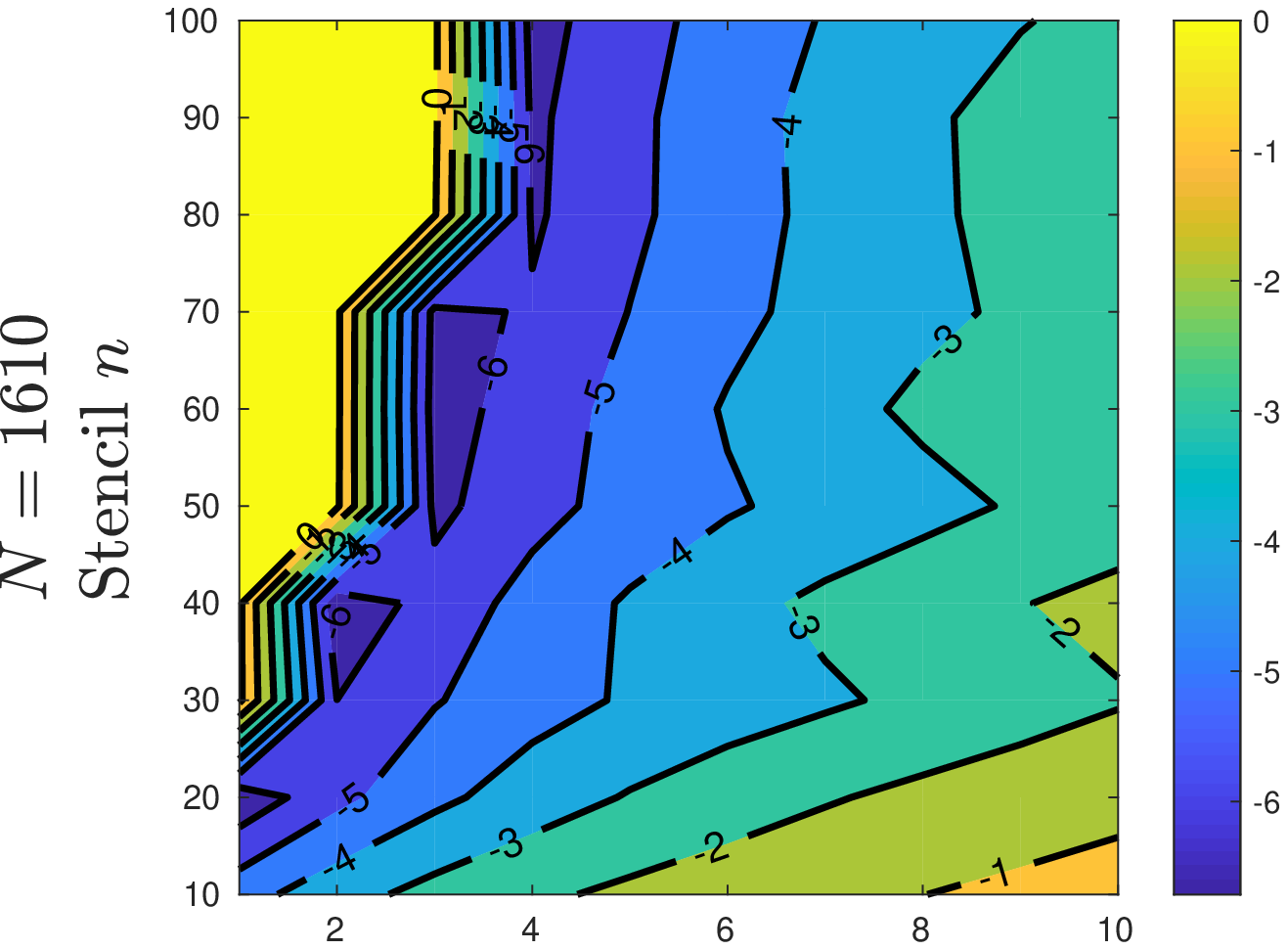}
\includegraphics[scale=.45]{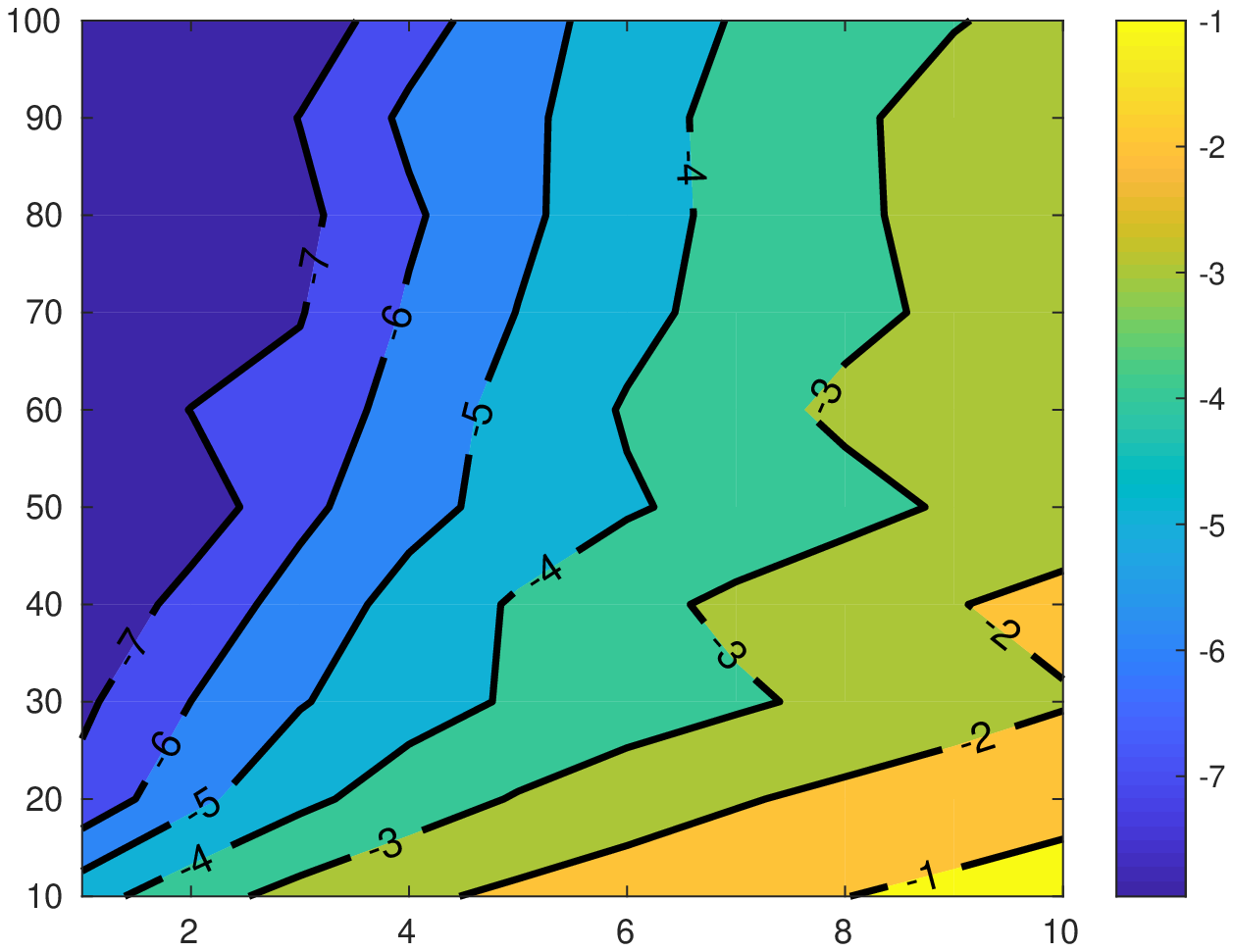}\\
\includegraphics[scale=.45]{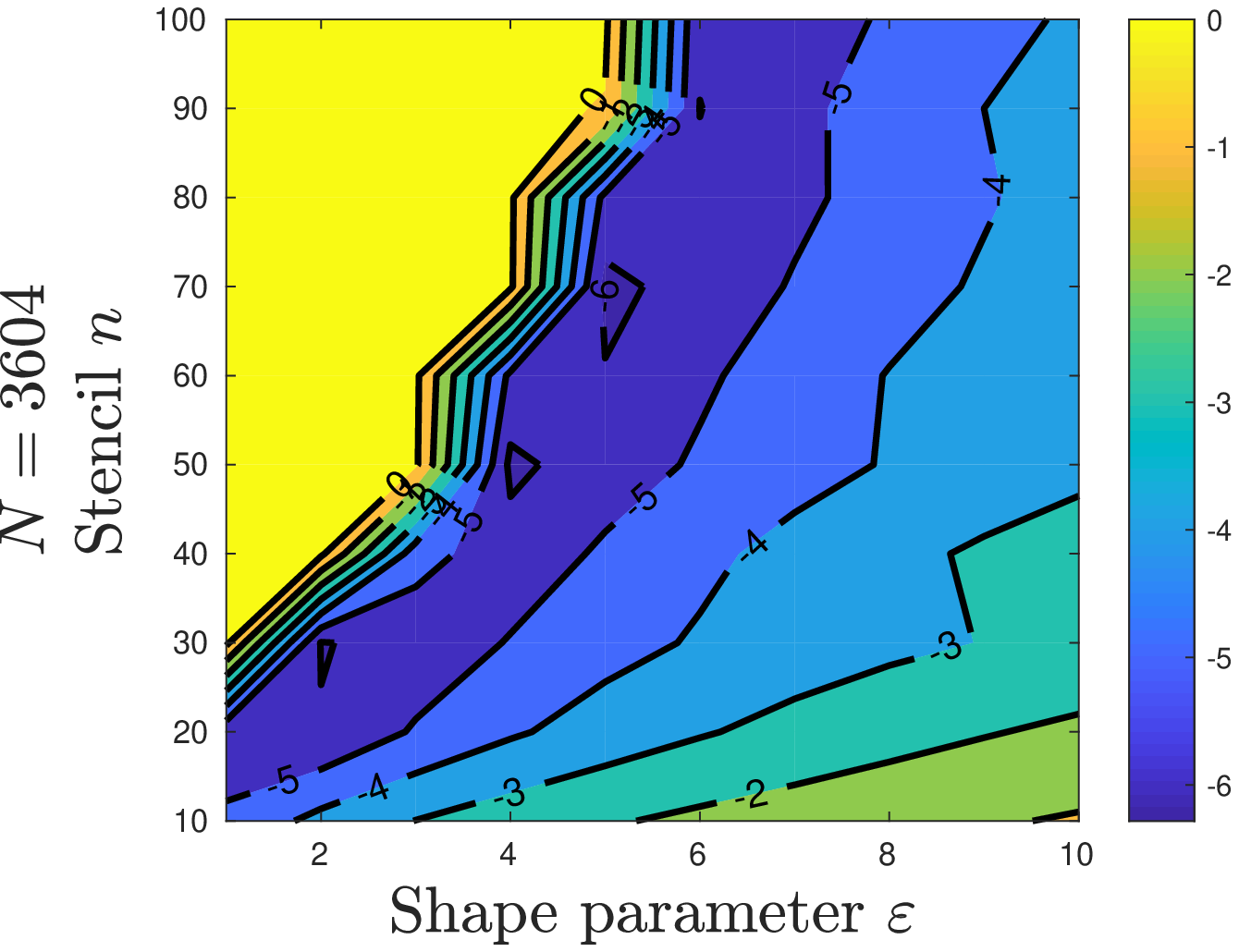}
\includegraphics[scale=.45]{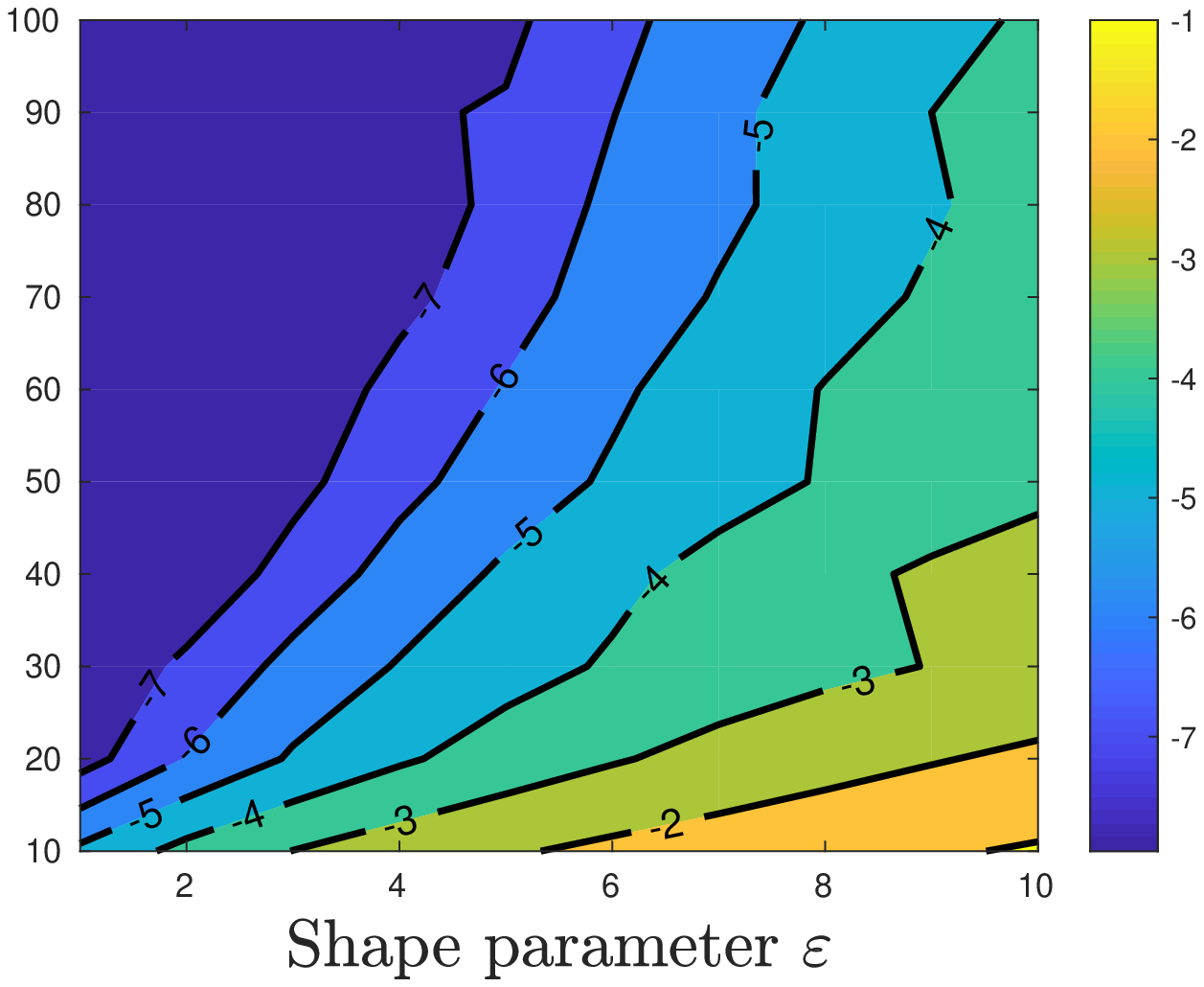} 
\caption{Accuracy isolines ($log_{10}$($L_2$-Error)) with $N_{int}=916,1610,3604$ interior points varying the shape parameter $\varepsilon$ and the stencil size $n$.}
\label{fg:HE_P7_isolines_quni916_error}
\end{figure*}

In Figure \ref{fg:HE_P7_isolines_quni916_error} we show the isolines of the error $log_{10}$($L_2$-Error) for the range of the shape parameter $[1,10]$ and for different sizes of stencils $n$=10:10:100. As $n$ increases, the linear systems increase, worsening the conditioning of the interpolation matrices. To understand the importance of local stability technique, both graphs of this figure must be observed simultaneously. The yellow region at the top left shows the region of error instability due to poor numerical conditioning while in the isolines of the graphs on the right, the region dark blue shows how $1\times 10^{-8}$ could be kept in order. As $N$ increases from 916 to 3604 this numerical behaviour is similar reading the figure row-wise.

\newpage

In \cite{mishra_fasshauer_sen_ling_2019} this same Helmholtz type PDE is worked with mixed type BC. In said work it can be seen that for $N=900$ nodes the $L_2$-Error $1\times 10^{-5}$ is reached using the Radial Basis Function - Finite Difference (RBF-FD) technique using a kernel hybrid of the Gaussian of type $\phi(r)=\alpha e^{-(\varepsilon r)^2}+\beta r^3$.

\newpage

\begin{figure*}[htb!!!!]
\includegraphics[scale=.45]{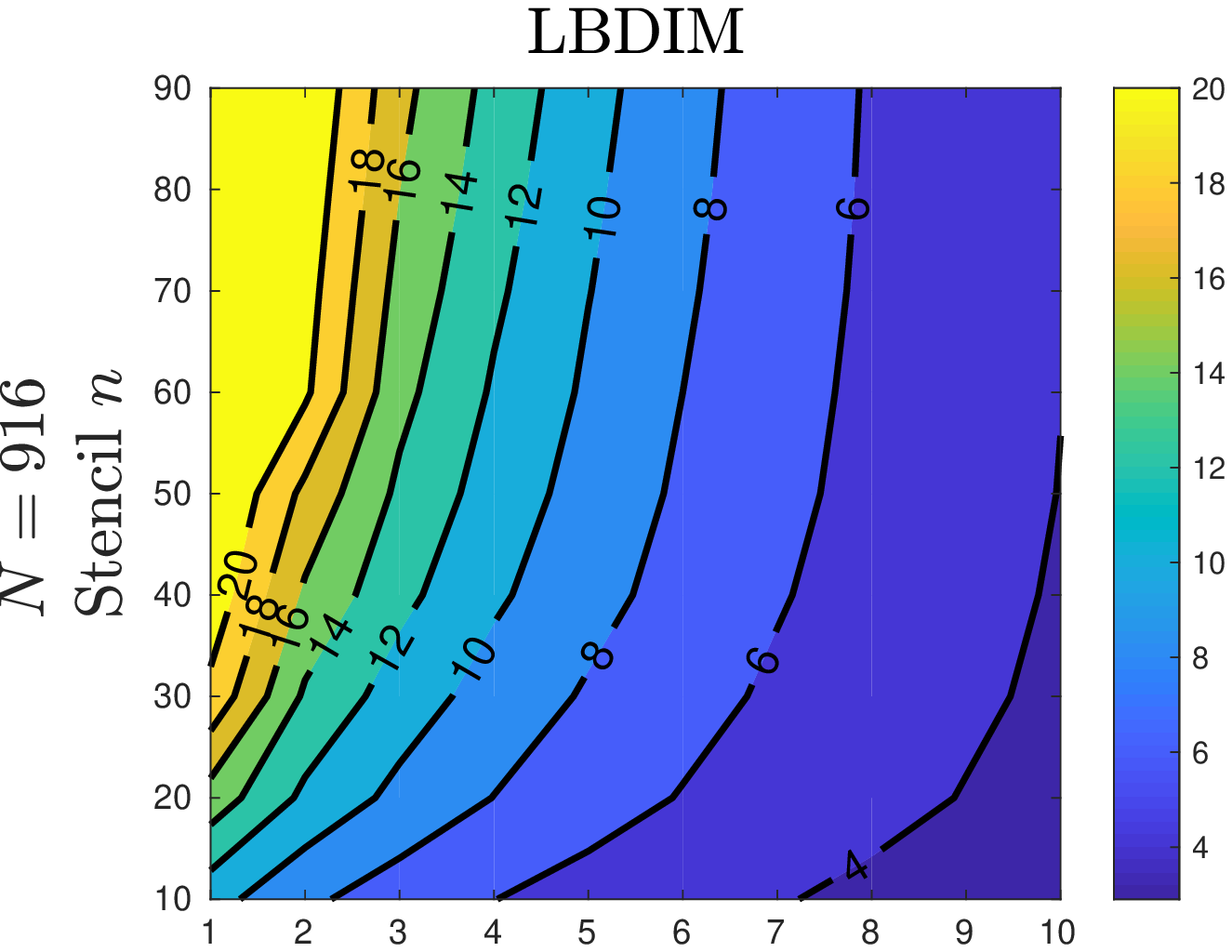}
\includegraphics[scale=.45]{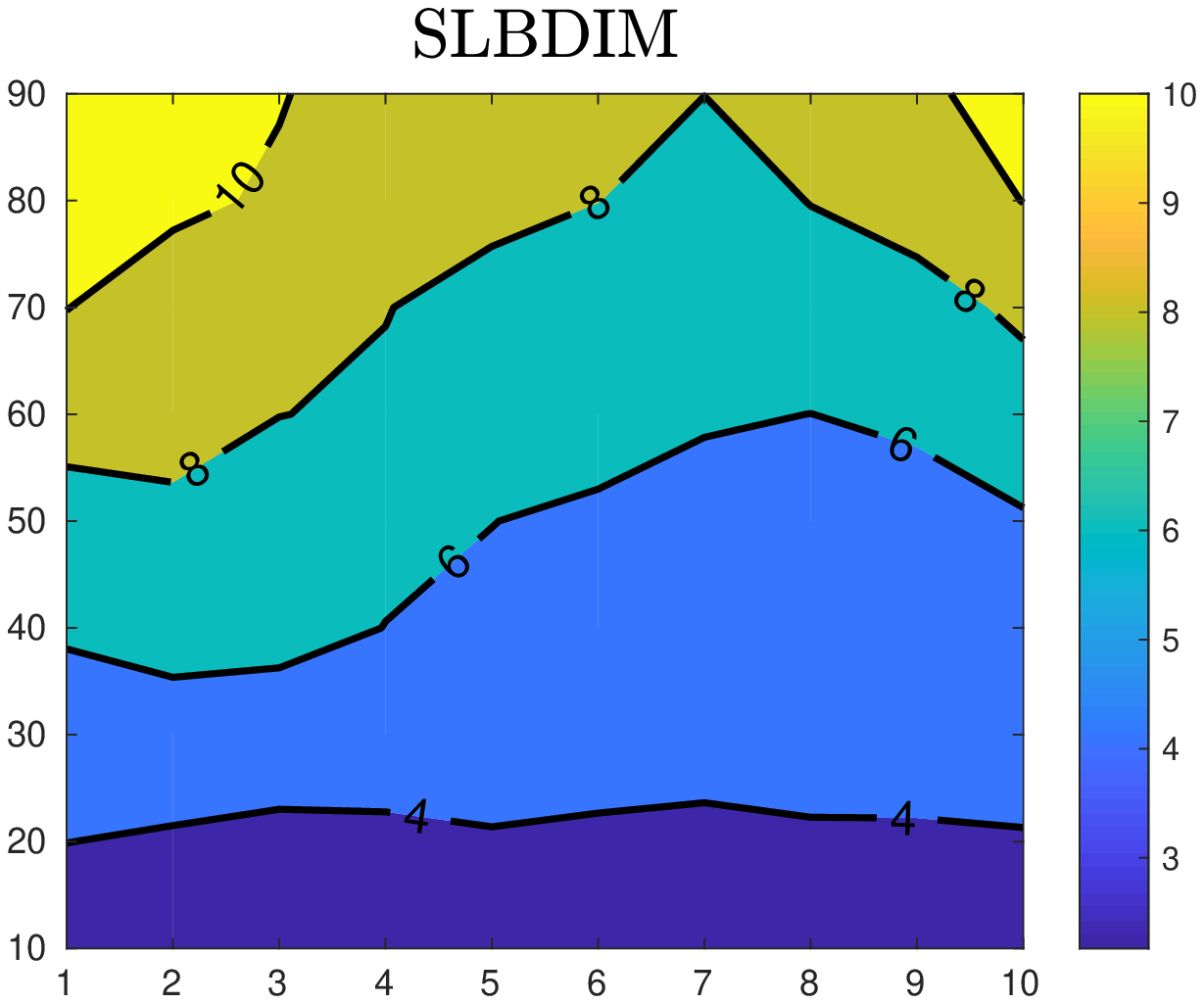}\\
\includegraphics[scale=.45]{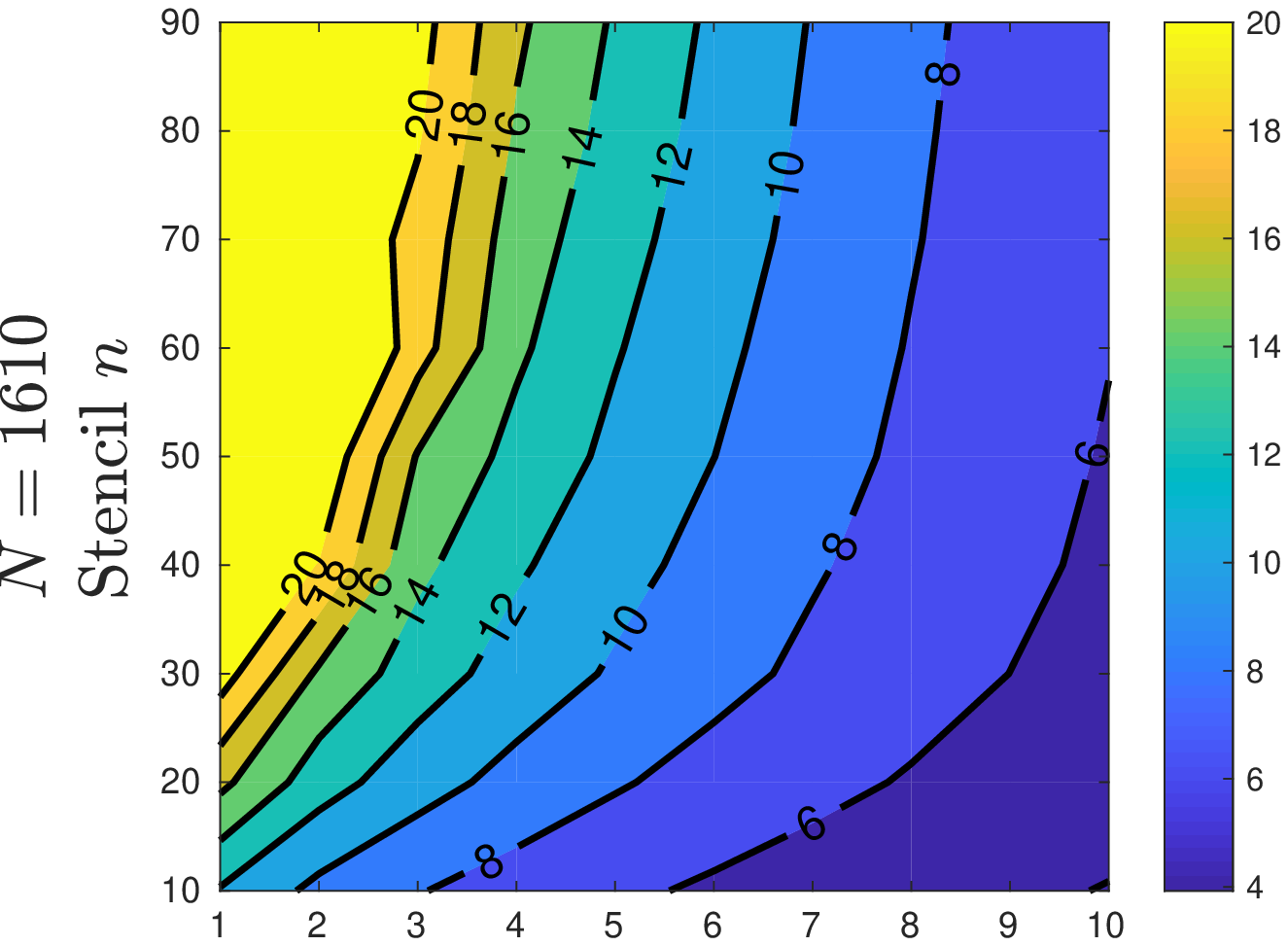}
\includegraphics[scale=.45]{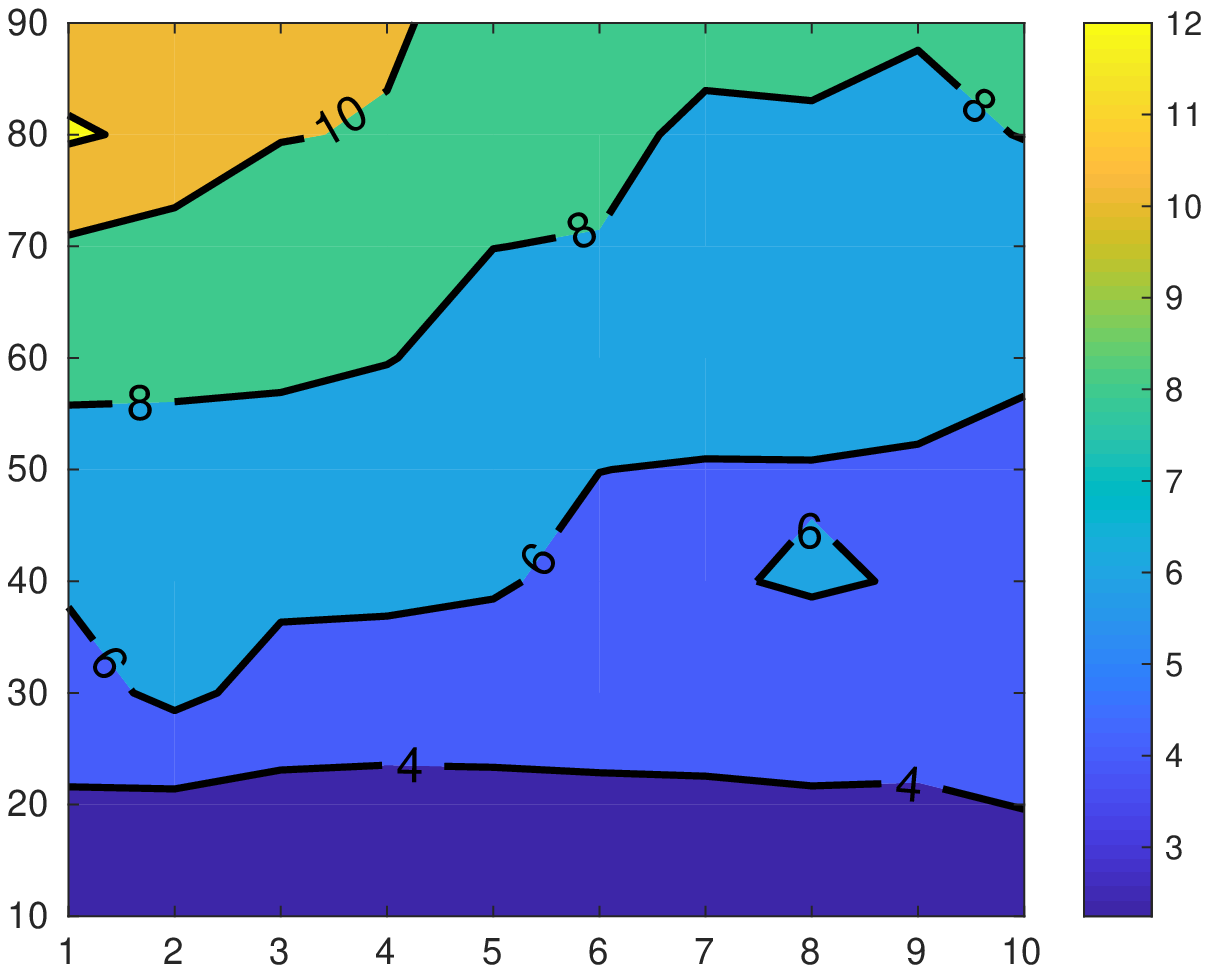}\\
\includegraphics[scale=.45]{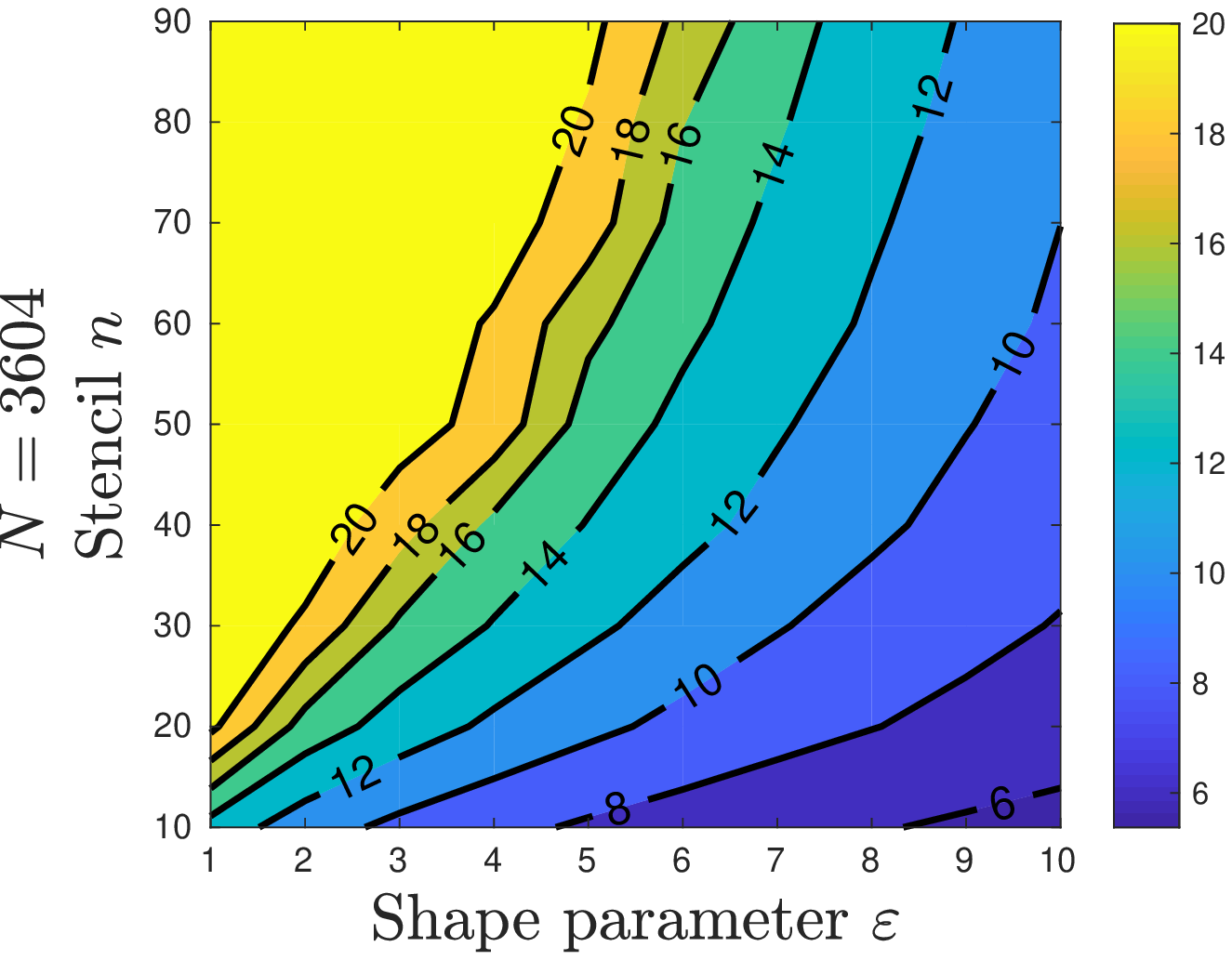}
\includegraphics[scale=.45]{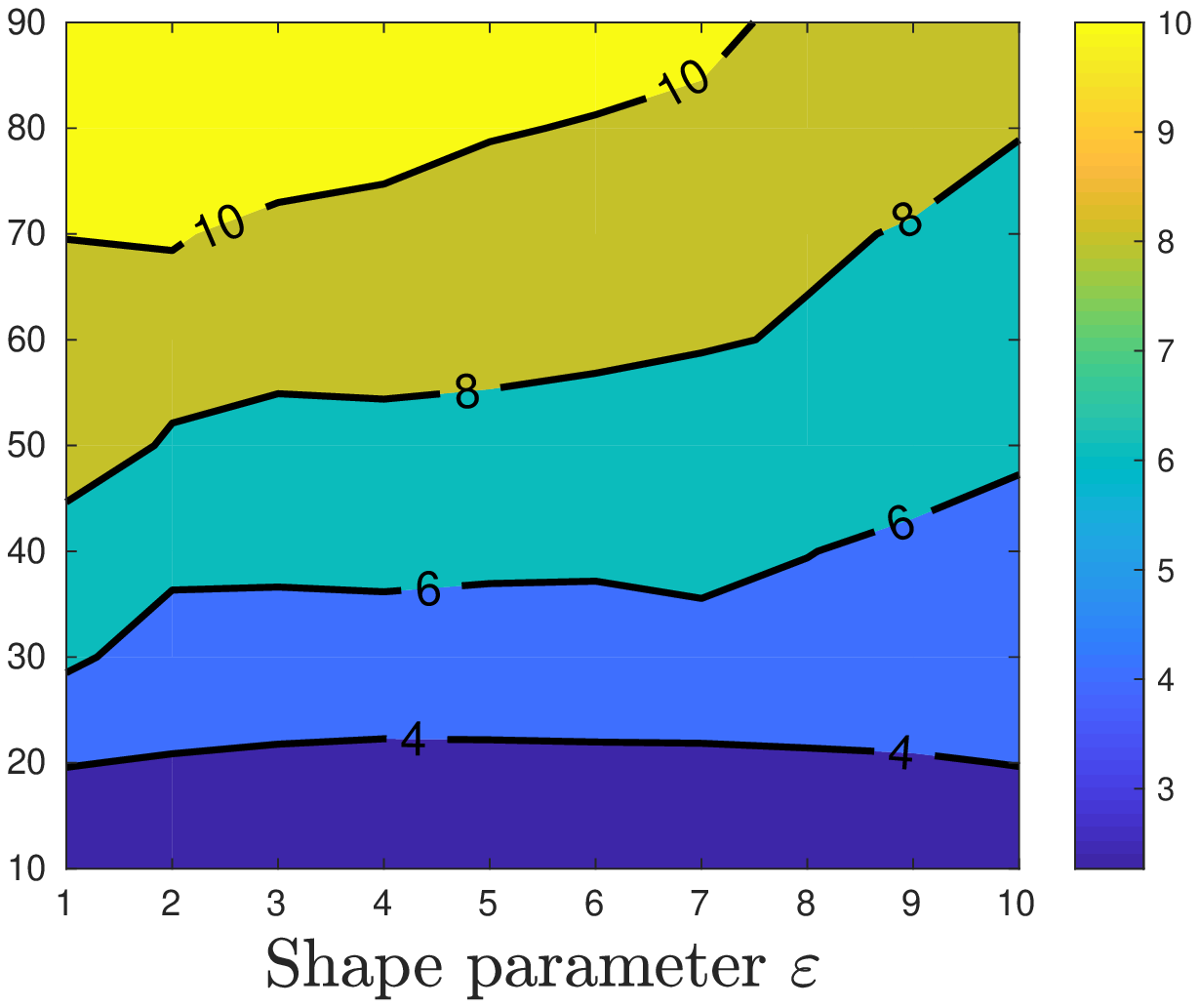} 
\caption{Condition number isolines ($log_{10}(\kappa(\boldsymbol{A}_i)$) with $N_{int}$=916,1610,3604 interior points varying the shape parameter $\varepsilon$ and the stencil size $n$.}
\label{fg:HE_P7_isolines_quni916_cn}
\end{figure*}

In Figure \ref{fg:HE_P7_isolines_quni916_cn} the isolines condition number $log_{10}(\kappa(\boldsymbol{A}_i)$ is shown. The range of the shape parameter is $[1,10]$ and the for different sizes of stencils are $n$=10:10:100. As $n$ increases, the conditioning of the local interpolation matrices increases. The yellow region at the top left shows the region of the condition number up to $1\times 10^{20}$. In the isolines of the graphs on the right column, the region dark blue shows better conditioning up to $1\times 10^{10}$. This ten order of magnitude are significant when when using linear solvers numerically. Also we can observe thar as $N$ increases from 916 to 3604 the conditioning behaviour is similar reading the figure row-wise.

In Figure \ref{fg:HE_P7_isolines_quni916_error} it was observed that the error plots suggest smaller values of $\varepsilon0$ for better accuracy, whereas in this figure the condition isolines plots suggest the larger values of $\varepsilon$ for better stability. This numerical results are interpreted as the well-known uncertainty principle in RBF local interpolations \cite{schaback_1995}. The idea behind this principle is that one cannot simultaneously achieve good conditioning and high accuracy using RBF basis. The relation between numerical stability and accuracy may be reviewed from different perspectives as in our case to obtain a stable formulation our option was to find a better basis in the same space of approximation using RBF-QR \cite{fornberg_larsson_flyer_2011} in the local boundary domain integral method.

\subsection{Polygonal billiars: case 2}
\label{subsec:polygo2}

\begin{figure*}[htb]
\centerline{\includegraphics[scale=.45]{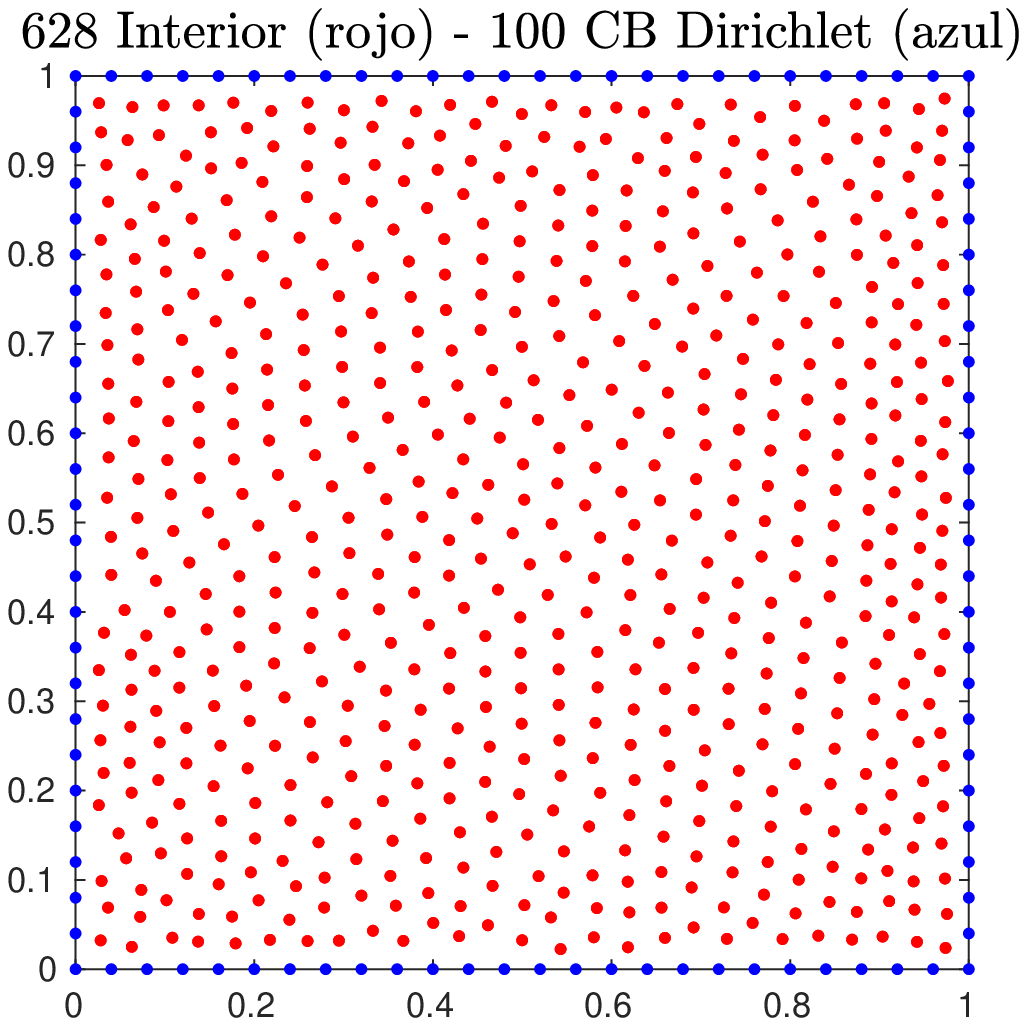}
\hspace{1cm}
\includegraphics[scale=.45]{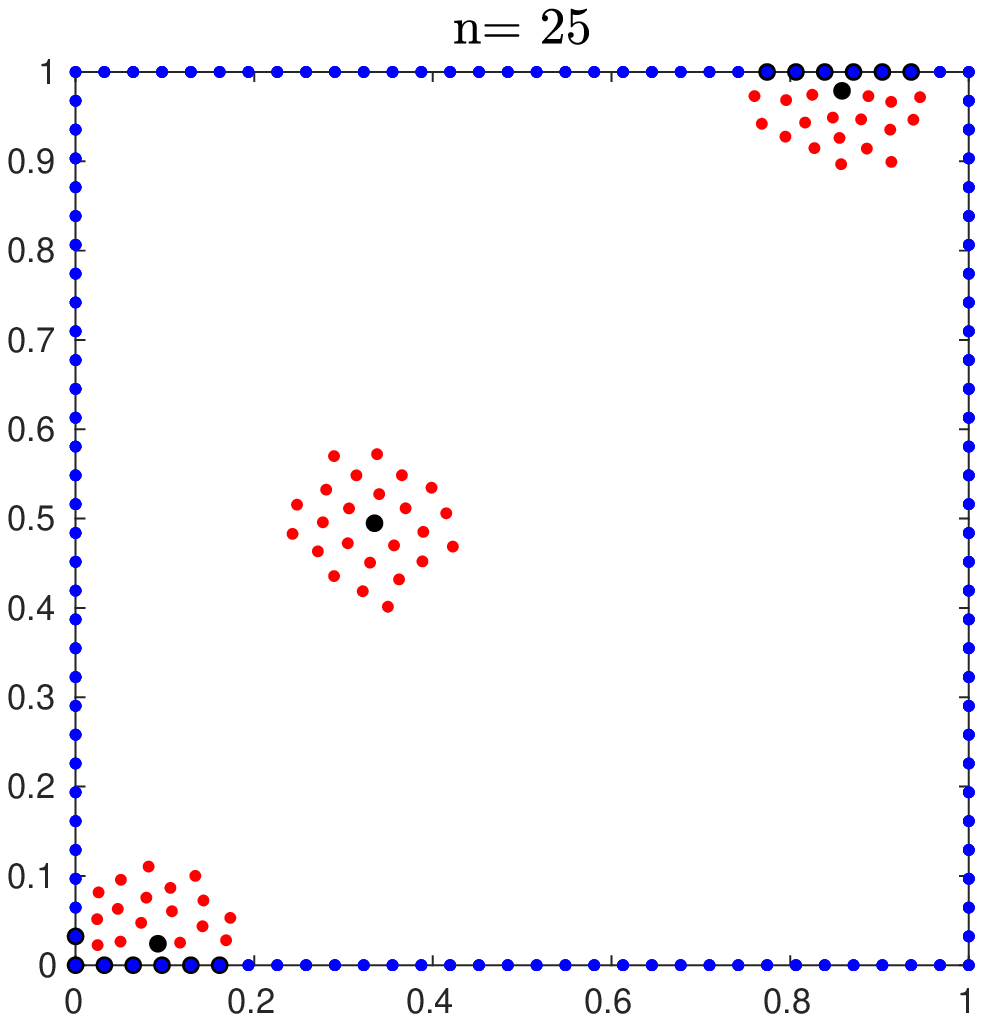}} \caption{Quasi-uniform node distribution with $N=900$ interior nodes (left). Stencil node sets with $n=25$ (right).} 
\label{fg:HE_P13_quni900_estencil25}
\end{figure*}

Consider the following two-dimensional Helmholtz equation
\begin{equation}
 \left\{
 \begin{array}{rcl}
  \Delta u(x,y) + k^2 u(x,y) &=& f(x,y), \hspace{.5cm}  \Omega =[0,1]\times [0,1],\\
  u(x,y) &=& g(x,y), \hspace{.5cm} \Gamma=\partial\Omega,
\label{Eq:edp_helmohltz_2}
 \end{array}
\right.
\end{equation}
where $k^2=2$, $f(x,y)=2x-4y$ and the exact solution is given by $u(x,y)=sin(\sqrt{3}x)sinh(y)+cos(\sqrt{2}y)+x-2y$, and $g (x,y)$ is chosen to match the exact one, thus giving BC of type Dirichlet. We use quasi-uniform nodes within the domain and stencils of size $n=25$ counting the collocation center as shown in Figure \ref{fg:HE_P13_quni900_estencil25}.

In Table \ref{tab:he_p13_ecm} we show the accuracy of the SLBDIM for the shape parameter $\varepsilon =1$ and for a range of low values, $\varepsilon\in \{ 1\times 10 ^0, 1\times 10^{-1}, 1\times 10^{-2}, 1\times 10^{-3}, 1\times 10^{-4}, 1\times 10^{- 5}\}$. The number of quasi-uniform interior points of the domain, $N$, varies from 121 to 900. It can be seen that for fixed $\varepsilon =1$, the \'orders of magnitude decrease from $1\times 10^{-6}$ to $1\times 10^{-8}$ starting at 441 nodes. In turn, the convergence of the method is observed for low values of the shape parameter, reaching RMS of the order $1\times 10^{-8}$ from 225 nodes. The $\varepsilon$ shown is where the best error is reached in that range.

\begin{table}[ht!]
\centering
\begin{tabular}{|c|c c|c c|} 
\hline
  $N$ & \multicolumn{2}{c|}{SLBDIM} & \multicolumn{2}{c|}{SLBDIM} \\
      &$\epsilon$ & $RMS$ & low $\epsilon$ & $RMS$\\
\hline
 $121$  & 1.0 & 1.2028E-06 & 0.1 & 2.1405E-07 \\ 
 $225$  & 1.0 & 5.8570E-07 & 0.1 & 5.0834E-08 \\ 
 $361$  & 1.0 & 3.9338E-07 & 0.01 & 3.3821E-08 \\ 
 $441$  & 1.0 & 7.8581e-08 & 0.1  & 3.3866E-08 \\
 $530$  & 1.0 & 5.2907E-08 & 0.00001 & 3.5984E-08 \\ 
 $628$  & 1.0 & 4.3843E-08 & 0.00001 & 3.6887E-08 \\ 
\hline
\end{tabular}
\caption{RMS for low shape parameters $\varepsilon\in \{1\times 10^{-1},\dots, 1\times 10^{-5}\}$.}
\label{tab:he_p13_ecm}
\end{table}

In \cite{ling_chen_sze_2012} this differential problem with mixed BC over the same domain is investigated using Multiquadric RBF kernels $\varphi(r,\varepsilon)=\sqrt{1+(\varepsilon r)^ 2}$ and a new RBF with $N\in [50,350]$ placement points. The results obtained in said reference reach errors of the order of $1\times 10^{-5}$ for $\varepsilon\in [0.4]$.





\section{Summary}
\label{sec:summ}

In this work we have introduced a new local integral method to compute resonances in dielectric cavities with different shapes. We have discussed numerical solutions, the node quasi-uniform node distributions over the domains and cavities with corners. Numerical results for Helmholtz-type equations were obtained using a stabilized local integral method that uses interpolations with RBF Gaussians. This method does not depend on a mesh, so it can be easily adapted to problems with complex geometries from . The good performance of the method has been shown with good results as shown in numerical tests 1 and 2 comparing with other results in the literature. Test 1 shows the advantage of using the SLBDIM to find regions of convergence of the $L_2$-Error of the order $1\times 10^{-8}$ when the shape parameter approaches zero. In test 2, a low shape parameter range is studied reaching the same order of the RMS. Having investigated the computational efficiency of the method, the future work consists of approaching some applications in wave chaos and dielectric microresonators, which is adequate to deal with geometries that come from arbitrary domains without analytical solutions.



\bibliographystyle{abbrv}
\bibliography{PoRa-article-references}
\end{document}